\documentclass[a4paper]{article}
\usepackage{amsmath,amssymb,amsthm,amsfonts,amscd}

\def\N{\mathbb{N}}
\def\gsum{\sum\limits}

\def\eps{\varepsilon}
\theoremstyle{plain}
\newtheorem{theorem}{Theorem}
\newtheorem{lemma}{Lemma}

\begin{document}
\author{Evgeny Shchepin}

 \title{ Greedy Sums and Dirichlet Series.}
\maketitle

\paragraph{Numerical arrays.} Under \emph{array of numbers} $\{a_i\}_{i\in
I}$ (over $I$) we mean a set of complex or real numbers indexed by a set  $I$ of arbitrary nature.
In particular, one can consider any  series as array over  natural numbers $\N$.

\emph{The direct product} of  arrays $\{a_i\}_{i\in
I}$ and $\{b_j\}_{j\in
J}$ is defined as array of products $\{a_ib_j\}_{i,j\in
I\times J}$.
In particular, the direct product of series forms a double series.

\paragraph{Greedy sum.}

For a numerical array $\{a_i\}_{i\in
I}$ and any positive $\eps$ we define its \emph{$\eps$-partial sum}
as $\sum\limits_{|a_i|\ge\eps}a_i $. Such sum exists iff it contains finite number of summands.

The limit of $\eps$-partial sums   $\lim_{\eps\to0}
\sum\limits_{|a_i|\ge\eps}a_i $, is called \emph{the greedy sum} of the array $\{a_i\}_{i\in
I}$.

The greedy sum of array $\{a_i\}_{i\in
I}$ is denoted by
 $$\gsum_{i\in I} a_i$$

The array $\{a_i\}_{i\in
I}$ having the greedy sum is called
\emph{greedy summable}.

 The greedy sum of any absolutely convergent series obviously coincides with its usual sum.
But the sum of a conditionally convergent series may differs from its greedy sum.

\paragraph{The main result.}
\begin{theorem}\label{t-product} If numeric arrays  $\{a_i\}_{i\in I}$, $\{ b_j\}_{j\in J}$ and their direct product
 $\{a_ib_j\}_{(i,j)\in I\times J}$ are greedy summable, then one has
\begin{equation}\label{product}
    \sum\limits_{(i,j)\in I\times J} a_ib_j=\sum\limits_{i\in I} a_i
    \sum\limits_{j\in J} b_j
\end{equation}
\end{theorem}

The proof of the theorem \ref{t-product} is based on the concept of \emph{greedy zeta-function} of array.
For a given numeric array $\{a_i\}_{i\in I}$ the function of complex variable $A(z)$ defined as
\begin{equation}\label{gen-funct}
  A(z)=\gsum_{i\in I}a_i|a_i|^z
\end{equation}
is called \emph{greedy zeta-function} of the array.
The theorem \ref{t-product} and the following equalities
\begin{equation}\label{ee}
  a_i|a_i|^z\cdot b_j|b_j|^z=a_ib_j|a_ib_j|^z
\end{equation}
 immediately implies the following theorem:
\begin{theorem}\label{z-product} Greedy zeta-function of the direct product of two arrays is equal to the product of
their greedy zeta-functions.
\end{theorem}
On the other hand the theorem \ref{t-product} represents a particular case of the theorem \ref{z-product}
($z=0$).

The proof of the theorem \ref{z-product} is based on the theory of generalized Dirichlet series developed by G. Hardy and M. Riesz
in the book \cite{H-R}.

\paragraph{The Dirichlet series of array.}

By \emph{Dirichlet series} we mean so called generalized Dirichlet series (see \cite{H-R})., i.e. a series of the form
\begin{equation}\label{Dirichlet}
\sum\limits_{n=1}^\infty c_n e^{-\lambda_n z},
\end{equation}
where  $\lambda_n$ --- \emph{exponents of series} represent a
monotone increasing to infinity sequence of real numbers and
\emph{coefficients of the series} $c_n$ as well as the variable
$z=x+iy$ are complex.

For a numeric array $\{a_i\}_{i\in I}$ let us define \emph{
Dirichlet series} $\sum c_k e^{-\ lambda_k z}$ in the following way:
$-\lambda_k$ is defined as  $k$-th by value element of the set
$\{\ln |a_i|\mid i\in I\}$. In particular, $-\lambda_1$ is the maximal element of the set
 $\{\ln |a_i|\mid i\in I\}$.
The coefficient  $c_k$ is defined as the sum
$\sum\limits_{\ln|a_i|=-\lambda_k} a_i$.

Immediately from the definition follows
\begin{lemma}
\label{greedy-dirichlet} Convergence of Dirichlet series of an
array $\{a_i\}_{i\in I}$ for a given value of variable $z$ is
equivalent  to existence of greedy sum $\gsum_{i\in I} a_i|a_i|^z$
and the greedy sum of the array is equal to the sum of
corresponding Dirichlet series for $z=0$.
\end{lemma}
Therefore greedy zeta-function of array coincides with its
Dirichlet series.

\paragraph{Multiplication theorem.}

One defines the \emph{formal product} of Dirichlet series $\sum
a_k e^{-\lambda_k z}$ and  $\sum b_k e^{-\mu_k z}$ as a Dirichlet
series $\sum c_k e^{-\nu_k z}$ such that $\nu_k=\lambda_i+\mu_j$
for some $i,j$ and $c_k=\sum\limits_{\lambda_i+\mu_j=\nu_k}
a_ib_j$.

It follows immediately from the definitions, that the Dirichlet
series of the direct product of two arrays is equal to the formal
product of Dirichlet series of the factors.

The theorem 55 of the book \cite{H-R} may be formulated as
follows.

\begin{theorem} Let us given two convergent series $\sum a_k$, $\sum b_k$
and two monotone increasing to infinity sequences $\lambda_k$ and
$\mu_k$ of positive numbers. Let $\sum c_k e^{-\nu_kz}$ the formal
product of Dirichlet series  $\sum a_k e^{-\lambda_kz}$ and $\sum
b_k e^{\mu_kz}$. If the series $\sum c_k$ converges then its sum
is equal to the product  $\sum a_k\sum b_k$.
\end{theorem}

To derive the theorem \ref{t-product} from the above theorem let
us simply apply it to the case $\lambda_k=-\ln |a_k|$, $\mu_k=-\ln
|b_k|$. In this case the corresponding Dirichlet series coincide
with corresponding greedy sums and its formal product corresponds
to greedy sum of the direct product of arrays,

\end{document}